\documentclass[12pt]{amsart}
\usepackage{amssymb, color}
\usepackage{amsmath}
\usepackage{amsthm}
 \font\tenmsb=msbm10 at 12pt \font\sevenmsb=msbm7 at 8pt \font\fivemsb=msbm5 at
6pt
\newfam\msbfam
\textfont\msbfam=\tenmsb \scriptfont\msbfam=\sevenmsb \scriptscriptfont\msbfam=\fivemsb
\def\Bbb#1{{\tenmsb\fam\msbfam#1}}

\def\Z{\Bbb Z}

\def\S{\Bbb S}

\begin{document}
\def \theequation{\arabic{section}.\arabic{equation}}
\newcommand{\reset}{\setcounter{equation}{0}}

\newcommand{\beq}{\begin{equation}}
\newcommand{\noi}{\noindent}
\newcommand{\eeq}{\end{equation}}
\newcommand{\dis}{\displaystyle}
\newcommand{\mint}{-\!\!\!\!\!\!\int}

\def \bx{\hspace{2.5mm}\rule{2.5mm}{2.5mm}} \def \vs{\vspace*{0.2cm}} \def
\hs{\hspace*{0.6cm}}
\def \ds{\displaystyle}
\def \p{\partial}
\def \O{\Omega}
\def \b{\beta}
\def \m{\mu}
\def \l{\lambda}
\def \le{\lambda^*}
\def \ul{u_\lambda}
\def \D{\Delta}
\def \d{\delta}
\def \s{\sigma}
\def \e{\varepsilon}
\def \a{\alpha}
\def \g{\gamma}
\def \R{\mathbb{R}}
\def \S{\mathbb{S}}
\def\qed{%
\mbox{ }%
\nolinebreak%
\hfill%
\rule{2mm} {2mm}%
\medbreak%
\par%
}
\newtheorem{thm}{Theorem}[section]
\newtheorem{lem}[thm]{Lemma}
\newtheorem{cor}[thm]{Corollary}
\newtheorem{prop}[thm]{Proposition}
\theoremstyle{definition}
\newtheorem{defn}{Definition}
\newtheorem{rem}[thm]{Remark}
\newenvironment{thmskip}{\begin{thm}\hfill}{\end{thm}}
\def \pr {\noindent {\it Proof:} }
\def \rmk {\noindent {\it Remark} }
\def \esp {\hspace{4mm}}
\def \dsp {\hspace{2mm}}
\def \ssp {\hspace{1mm}}

\title{ Existence of stable solutions to $(-\Delta)^m u=e^u$ in $\R^N$ with $m \geq 3$ and $N > 2m$}
\author{Xia Huang}
\address{Department of Mathematics and Center for Partial Differential Equations, East China Normal
University, Shanghai 200062, P.R. China} \email
{xhuang1209@gmail.com}
\author{Dong Ye}
\address{IECL, UMR 7502, D\'epartement de Math\'ematiques, Universit\'e de Lorraine, B\^at. A, \^{\i}le de Saulcy, 57045 Metz, France} \email{dong.ye@univ-lorraine.fr}
 \subjclass{35J91, 35B08, 35B35} \keywords{polyharmonic equation, entire stable solution}

\date{}
\begin{abstract}
We consider the polyharmonic equation $(-\Delta)^m u=e^u$ in $\R^N$ with $m \geq 3$ and $N > 2m$. We prove the existence of many entire stable solutions. This answer some questions raised by Farina and Ferrero in \cite{FF}.
\end{abstract}

\maketitle
\section{Introduction}

In this paper, we are interested in the existence of entire stable solutions of the polyharmonic equation
\begin{equation}\label{1.1}
(-\Delta)^m u=e^u\quad\text{in}~\mathbb{R}^N.
\end{equation}
with $m\geq 3$ and $N > 2m$.
\begin{defn}
{\it A solution $u$ to \eqref{1.1} is said stable in $\Omega\subseteq\mathbb{R}^N$ if
$$
\left\{\begin{array}{ll}
\displaystyle\int_\Omega|\nabla(\Delta^{\frac{m-1}{2}}\phi)|^2 dx-\int_\Omega e^u\phi^2 dx\geq 0\quad\text{for~any}~\phi\in C_0^\infty(\Omega),& \vspace*{0.2cm} \mbox{when $m$ is odd};\\
\displaystyle\int_\Omega|\Delta^{\frac{m}{2}}\phi|^2 dx-\int_\Omega e^u\phi^2 dx\geq0\quad\text{for~any}~\phi\in C_0^\infty(\Omega), & \mbox{when $m$ is even}.
\end{array}\right.
$$
Moreover, a solution to \eqref{1.1} is said stable outside a compact set $K$ if it's stable in
$\mathbb{R}^N\setminus{K}$. For simplicity, we say also that $u$ is stable if $\O = \mathbb{R}^N$.}
\end{defn}

For $m=1$, Farina \cite{F} showed that \eqref{1.1} has no stable classical solution in $\mathbb{R}^N$ for $1\leq N\leq 9$. He also proved that any classical solution which is stable outside a compact set in $\R^2$ verifies $e^u \in L^1(\R^2)$, therefore $u$ is provided by the stereographic projection thanks to Chen-Li's classification result in \cite{CL}, that is, there exist $\lambda>0$ and $x_0\in\mathbb{R}^2$ such that
\begin{align}
\label{new1.3}
u(x)=\ln\left[\frac{32\lambda^2}{(4+\lambda^2|x-x_0|^2)^2}\right]\quad\text{for~some}~\lambda>0.
\end{align}
Later on, Dancer and Farina \cite{DF} showed that \eqref{1.1} admits classical entire solutions which are stable outside a compact set of $\mathbb{R}^N$ if and only if $N\geq 10$.

\medskip
It is well known that for any $m \geq 1$, $\l > 0$ and $x_0 \in \R^{2m}$, the function $u$ defined
in \eqref{new1.3} resolves \eqref{1.1} in the conformal dimension $\R^{2m}$, they are the so-called spherical
solutions, since they are provided by the stereographic projections.

\medskip
 For $m=2$, the stability properties of entire solutions to \eqref{1.1} were studied in many works, especially the study for radial solutions is complete. Let $u(x)=u(r)$ be a smooth radial solution to \eqref{1.1}, then $u$ satisfies the following initial value problem
\begin{equation}\label{1.2}
\begin{cases}
\begin{aligned}
& (-\Delta)^m u=e^u, && \\
& u^{(2k+1)}(0)=0, &&  \forall\; 0\leq k\leq m-1,\\
& \Delta^ku(0)=a_k,&& \forall\; 0\leq k\leq m-1.
\end{aligned}
\end{cases}
\end{equation}
Here the Laplacian $\Delta$ is seen as $\Delta u = r^{1-N}\left(r^{N-1}u'\right)'$
and $a_k$ are constants in $\R$. Equivalently, let $v_k = (-\Delta)^ku$ for $0\leq k \leq m-1$, the equation \eqref{1.2} can be written as a system
 \begin{align}
 \label{new1.2}
-v_k'' - \frac{N-1}{r}v_k' = v_{k+1} \;\mbox{ for } \; 0\leq k \leq m-2;\quad \mbox{and} \;\;-v_{m-1}'' - \frac{N-1}{r}v_{m-1}' = e^{v_0}
 \end{align}
 where $v_k(0) = (-1)^ka_k$ and $v_k'(0) = 0$ for any $0\leq k \leq m-1$.

\medskip
Let $m = 2$, $a_0 = u(0) = 0$ (It's always possible by the scaling $u(\l x) + 2m\ln\l$).
Denote by $u_\beta$ the solution to \eqref{1.2} verifying $a_1 = \beta$, it's known from \cite{BF, DG, W} that:
\begin{itemize}
\item There is no global solutions to \eqref{1.2} if $N \leq 2$.
\item For $N \geq 3$, there exists $\beta_0 < 0$ depending on $N$ such that the solution to \eqref{1.2} is globally defined, if and only if $\beta \leq \beta_0$.
\item If $N = 3$ or $4$, any entire solution $u_\beta$ is unstable in $\R^N$, but stable outside a compact set.
\item If $5\leq N\leq 12$, then $u_\beta$ is
  stable outside a compact set for every $\beta<\beta_0$ while $u_{\beta_0}$ is unstable outside every compact set.
\item If $5\leq N\leq 12$, there exists $\beta_1 < \beta_0$ such that $u_\beta$ is stable in $\R^N$, if and only if $\beta \leq \beta_1$.
\item If $N\geq 13$, $u_\beta$ is stable for every $\beta \leq\beta_0$.
\end{itemize}
Moreover, Dupaigne {\it et al.} showed in \cite{DG} the examples of non radial stable solutions for $\D^2 u = e^u$
in $\R^N$ with any $N \geq 5$, and Warnault proved in \cite{W} that no stable (radial or not) smooth solution exists for $\D^2 u = e^u$ if $N \leq 4$.

\medskip
Recently, Farina and Ferrero \cite{FF} studied \eqref{1.1} for general $m \geq 3$, they obtained many results about
the existence and stability of solutions, especially for the radial solutions. More precisely, they proved that
\begin{itemize}
\item For $N\leq 2m$, no stable solution (radial or not) exists;
\item For $m\geq 3$ odd, if $1\leq N\leq 2m-1$ or $m\geq 1$ odd and $N=1$, then any radial solution is stable outside a compact set;
\item For $m\geq 1$ and $N=2m$, then the spherical solutions, i.e.~solutions given by \eqref{new1.3} are stable outside a compact set.
\item For $m\geq 3$ odd, if $(-1)^ka_k \leq 0$ for same $1\leq k\leq m-1$, then the radial solution is stable outside a compact set;
\item For $m\geq 2$ even and $u(0)=0$, there exists a function $\Phi: \R^{m-1} \rightarrow (-\infty, 0)$ (depending on $N$) such that the solution to \eqref{1.2} is global if and only if $a_{m-1} \leq \Phi(a_1,...,a_{m-2})$. Moreover, if $a_{m-1} < \Phi(a_1,...,a_{m-2})$, then the solution is stable outside a compact set;
\end{itemize}

It is also worthy to mention that for the conformal or critical dimension $N = 2m$ with $m \geq 2$, many existence results exist by prescribing the behavior of $u$ at infinity. See \cite {CC, WY, DG} for $m = 2$ and see \cite{HM} for $m \geq 3$. Clearly, these results imply the existence of many non radial solutions which are stable outside a compact set.

\medskip
However, in the supercritical dimensions $N>2m$ with $m \geq 3$, less is known for the stable solutions. Farina and Ferrero raised then the question
(see for instance Problem 4.1 (iii) in \cite{FF}) about the existence of stable solutions. In this work, we will provide rich examples of stable solutions. First we consider radial solutions to \eqref{1.2} and show that the solution is stable if we allow $a_{m-1}$ to be negative enough.
\begin{thm}
\label{thm1.1}
Let $m\geq 2$ and $N>2m$. Given any $(a_k)_{0\leq k \leq m-2}$, there exists $\beta \in \R$ such that the solution to equation \eqref{1.2} is stable in $\R^N$ for any $a_{m-1} \leq \beta$.
\end{thm}

Furthermore, given any $N > 2m$, we prove the existence of non radial stable solution to \eqref{1.1} and the existence of stable radial solutions for the following {\it borderline} situations, see Theorem \ref{thm1.2} and Corollary \ref{cor1.6} below.
\begin{itemize}
\item[(i)] $N > 2m$, $m \geq 3$ is odd, and $(-1)^ka_k > 0$ for any $1\leq k\leq m-1$;
\item[(ii)] $N > 2m$, $m \geq 4$ is even, $u(0) = 0$ and $a_{m-1} = \Phi(a_1,...,a_{m-2})$;
\end{itemize}
The existence of stable radial solutions on the borderline for $m\geq 4$ even in arbitrary supercritical dimension is a new phenomenon comparing to $m=2$, where the borderline solutions are not stable out of any compact set if $5\leq N \leq 12$.

\medskip
It will be interesting to know if all radial solutions are stable in high dimensions as for $m=2$ and $N \geq 13$. We are
not able to answer this question, but we can prove that for $m \geq 3$ odd, and a wide class of initial
data $(a_k)$, the corresponding radial solutions are effectively always stable in large dimensions.
\begin{thm}
\label{thm1.4}
Let $m\geq 3$ be odd, then there exists $N_0$ depending only on $m$ such that for any $N\geq N_0$, the radial solution to \eqref{1.2} with $a_k\leq 0$ for $1\leq k\leq m-1$ is stable in $\mathbb{R}^N$.
\end{thm}

The following Hardy inequalities will play an important role in our study of stability, see Theorem 3.3 in \cite{M}. Let $m\geq2$ and $N > 2m$. If $m$ is odd, then
\begin{equation*}
\l_{N,m}\int_{\mathbb{R}^N}\frac{\varphi^2}{|x|^{2m}} dx\leq\int_{\mathbb{R}^N}|\nabla(\Delta^{\frac{m-1}{2}}\varphi)|^2 dx\quad\text{for~any}~\varphi\in C_0^\infty(\mathbb{R}^N),
\end{equation*}
where
\begin{equation}\label{1.3}
\l_{N,m}:=\frac{(N-2)^2}{16^\frac{m}{2}}\prod_{i=1}^{\frac{m-1}{2}}(N-4i-2)^2(N+4i-2)^2.
\end{equation}
If $m$ is even, then
\begin{equation*}
\mu_{N,m}\int_{\mathbb{R}^N}\frac{\varphi^2}{|x|^{2m}} dx\leq\int_{\mathbb{R}^N}|\Delta^{\frac{m}{2}}\varphi|^2 dx\quad\text{for~any}~\varphi\in C_0^\infty(\mathbb{R}^N),
\end{equation*}
where
\begin{equation}\label{1.4}
\mu_{N,m}:=\frac{1}{16^{\frac{m}{2}}}\prod_{i=0}^{\frac{m-2}{2}}(N+4i)^2(N-4i-4)^2.
\end{equation}

We will also use the following well-known comparison result (see for instance Proposition 13.2 in \cite{FF})
\begin{lem}\label{lem1.2}
\label{comp}
Let $u, v\in C^{2m}([0, R))$ be two radial functions such that $\Delta^m u - e^u \geq \Delta^m v - e^v$ in $[0, R)$ and
\begin{align}
\Delta^ku(0) \geq \Delta^kv(0), \; (\Delta^ku)'(0) \geq (\Delta^kv)'(0), \quad \forall \; 0\leq k \leq m-1.
\end{align}
Then for any $r\in[0, R)$ we have
$$
\Delta^ku(r) \geq \Delta^kv(r),\quad \text{for~all}~0\leq k\leq m-1.
$$
\end{lem}

\section{A first existence result}
\reset
Here we prove Theorem \ref{thm1.1}. We will consider radial solutions to the initial value problem \eqref{1.2}. Denote
\begin{equation}\label{new1}
c_k=\Delta^k (r^{2k})=\prod_{i=1}^k 2i(N-2+2i)\quad\text{for~ any}~ k\geq 1.
\end{equation}

{\it Case 1:  $m\geq 3$ is odd.}

\smallskip
Fix $\Delta^ku(0)=a_k$ for $0\leq k\leq m-2$. Consider the solution $u_{(a_k)}$ to \eqref{1.2} associated to the initial values $a_k$, $0\leq k \leq m-1$. We know that the solution is globally defined in $\R^N$ for any $(a_k)$, see \cite{FF}. Clearly, the polynomial
$$\Psi(r) = a_0 + \sum_{1\leq k \leq m-1} \frac{a_k}{c_k}r^{2k} \quad \mbox{with $c_k$ given by \eqref{new1}}$$
verifies $\Delta^m \Psi \equiv 0$ in $\R^N$ and $\Delta^k\Psi(0)=a_k$ for all $0\leq k\leq m-1$.

\medskip
As $\Delta^m(u_{(a_k)}-\Psi) = -e^{u_{(a_k)}} < 0$, it's easy to check that $u_{(a_k)}(r) < \Psi(r)$ for any $r>0$. We claim that: $u_{(a_k)}$ is stable when $a_{m-1}$ is small enough. In fact, we need only to get the following estimate:
\begin{equation}\label{1.5}
e^{\Psi(r)}\leq \frac{\l_{N,m}}{r^{2m}}\quad \text{in}~\mathbb{R}^N,
\end{equation}
where $\l_{N,m}>0$ is given by \eqref{1.3}. Let
$$
h(r)={c_{m-1}}r^{2-2m}\left[a_0 + \sum_{1\leq k \leq m-2} \frac{a_k}{c_k}r^{2k}+ 2m\ln r - \ln \l_{N,m}\right].
$$
Obviously $\lim_{r\rightarrow+\infty} h(r)=0$ and $\lim_{r\rightarrow 0}h(r)=-\infty$. So $H_0 = \sup_{(0,\infty)} h(r) < \infty$ exists and \eqref{1.5} holds if $-a_{m-1} \geq H_0$. We conclude that if $a_{m-1} \leq -H_0$,
\begin{align*}
\int_{\mathbb{R}^N}|\nabla(\Delta^{\frac{m-1}{2}}\phi)|^2 dx-\int_{\mathbb{R}^N}e^{u_{(a_k)}}\phi^2 dx & \geq \int_{\mathbb{R}^N}|\nabla(\Delta^{\frac{m-1}{2}}\phi)|^2 dx-\int_{\mathbb{R}^N}e^{\Psi}\phi^2 dx\\
& \geq \int_{\mathbb{R}^N}|\nabla(\Delta^{\frac{m-1}{2}}\phi)|^2 dx-\l_{N,m}\int_{\mathbb{R}^N}\frac{\phi^2}{|x|^{2m}} dx\geq 0,
\end{align*}
i.e.~$u_{(a_k)}$ is stable in $\R^N$.

\medskip
{\it Case 2:  $m$ is even.}

\smallskip
Let $\Delta^ku(0)= a_k$ for $0\leq k\leq m-2$ be fixed. We can check that the scaling $u(\l x) + 2m\ln\l$ does not affect the stability of the solution, so we can assume that $a_0 = 0$ without loss of generality. By Theorem 2.2 in \cite{FF}, the solution to \eqref{new1.2} is global if and only if $a_{m-1} \leq \beta_0 = \Phi(a_k)$. For any $a_{m-1} < \beta_0$, consider
 $$
\Psi(r)=u_{\beta_0}(r)+\frac{(a_{m-1}-\beta_0)r^{2m-2}}{c_{m-1}},
$$
then $\Delta^m\Psi =\Delta^m u_{\beta_0} = e^{u_{\beta_0}} \geq e^\Psi$. Using Lemma \ref{lem1.2}, we have $u_{(a_k)} \leq \Psi$ in $\mathbb{R}^N$
as $\Delta^k\Psi(0)=\Delta^k u_{(a_k)}(0)$ for any $0\leq k\leq{m-1}$. As above, if there holds
\begin{equation}\label{1.7}
e^{u_{\beta_0}} e^{\frac{(a_{m-1}-\beta_0)r^{2m-2}}{c_{m-1}}}\leq\frac{\mu_{N,m}}{r^{2m}}\quad \mbox{in }\;\R^N,
\end{equation}
with $\mu_{N,m}$ given by \eqref{1.4}, then $u_{(a_k)}$ is stable in $\R^N$. Let $$
g(r)={c_{m-1}}r^{2-2m}\left[u_{\beta_0}(r)-\ln\frac{\mu_{N,m}}{r^{2m}}\right]-\beta_0.
$$
By \cite{FF}, the borderline entire solution $u_{\beta_0}(r) = o(r^{2m-2})$ as $r \to \infty$.
So $\lim_{r\rightarrow+\infty} g(r)= -\beta_0$, $\lim_{r\rightarrow0}g(r) =-\infty$, and \eqref{1.7} holds if we take $-a_{m-1} \geq\sup_{(0,\infty)} g$. \qed

\section{More stable solutions}
Here we will prove for $N > 2m$ the existence of stable radial solution to \eqref{1.2} for the borderline cases
but also the existence of non radial stable solutions to \eqref{1.1}. We begin with the case with $m$ odd.
\begin{thm}
\label{thm1.2}
For $m\geq 3$ be odd and $N>2m$, then there exists entire stable solution $u$ of \eqref{1.2} satisfying sign$(a_k)=(-1)^k$ for all $0 \leq k\leq m-1$.
\end{thm}

\noindent
Proof. Consider $u^\e$, solution of \eqref{1.2} with the initial conditions $a_k = (-1)^k\e$ for $0 \leq k \leq m-3$; $a_{m-2} = -\beta$ with $\beta > 0$ and $a_{m-1} = \e$. Here $\e \in (0, 1]$ is a small parameter, for simplicity, we will omit the exponent $\e$ in the following. Let
$$\Psi(r) := -\frac{\beta}{c_{m-2}}r^{2m-4} + \e H(r),
 $$
where
$$H(r) := 1 + \sum_{k=1}^{m-3} \frac{(-1)^k}{c_k}r^{2k} + \frac{r^{2m-2}}{c_{m-1}} \quad
\mbox{with $c_k$ given by \eqref{new1}.}$$
Therefore $(-\Delta)^m\Psi \equiv 0$ and $\Delta^k\Psi(0) = \Delta^k u(0)$ for any $0\leq k \leq m-1$.
Denote also $$H_+(r) := 1 + \sum_{k=1}^{m-3} \frac{r^{2k}}{c_k} + \frac{r^{2m-2}}{c_{m-1}}.$$
As we have
$$u \leq \Psi \leq -\frac{\beta}{c_{m-2}}r^{2m-4} + \e H_+(r) \quad \mbox{in }\; [0, \infty),$$ there holds
$u(r) \leq \e H_+(1)$ in $[0, 1]$. Denote $\gamma_0 := e^{H_+(1)}$ and consider
$v := u - \Psi + \frac{\gamma_0}{c_m}r^{2m}$. Then
$\Delta^m v = \Delta^m u + \gamma_0 = -e^u + \gamma_0 \geq 0$ for any $\e \leq 1$ and $r \in [0, 1]$.
Since $\Delta^k v(0) = 0$ for any $0\leq k \leq m-1$, we get $v \geq 0$ in $[0, 1]$, hence
\begin{align*}
u(r) & \geq \e H(r) - \frac{\beta}{c_{m-2}}r^{2m-4} - \frac{\gamma_0}{c_m}r^{2m} \\
& > - H_+(1) - \frac{\beta}{c_{m-2}} - \frac{\gamma_0}{c_m} =:\xi_0, \quad \forall\; r \in [0, 1], \;\e \leq 1.
\end{align*}
Inversely, consider $w := u - \Psi + \frac{e^{\xi_0}}{c_m}r^{2m}$ in $[0, 1]$, there holds $\D^m w = e^{\xi_0} - e^u \leq 0$ in $[0, 1]$. By the lemma \ref{1.2}, we have then $\Delta^k w(r) \leq 0$ in $[0, 1]$ for any $0 \leq k \leq m$, so that for $r\in [0, 1]$,
\begin{align*}
\Delta^{m-1}u(r) \leq \e - e^{\xi_0}\frac{r^2}{2N}, \quad \Delta^{m-2} u(r) \leq -\beta + \e\frac{r^2}{2N} - e^{\xi_0}\frac{r^4}{8N(N+2)}.
\end{align*}
Moreover, as $\Delta^{m-1} u$ is decreasing, we have $\Delta^{m-1}u(r) \leq \Delta^{m-1}u(1) \leq \e -\frac{e^{\xi_0}}{2N}$ in $(1, \infty)$. Consequently, for $r> 1$,
\begin{align*}
\Delta^{m-2} u(r) &=\Delta^{m-2} u(1) +\int_{1}^r\rho^{1-N}\int_0^\rho s^{N-1}\Delta^{m-1} u(s)  ds d\rho\\&\leq -\beta + \frac{\e}{2N} - \frac{e^{\xi_0}}{8N(N+2)} + \int_1^r\rho^{1-N}\int_0^\rho\left[\e - e^{\xi_0} \frac{\min(1, s)^2}{2N} \right]s^{N-1} ds d\rho\\& = -\beta + \e\frac{r^2}{2N} - e^{\xi_0}\left[ \frac{1}{8N(N+2)} + \frac{1}{2N^2}\int_1^r\left(\rho - \frac{2}{N+2}\rho^{1-N}\right) d\rho\right]\\
& = -\beta + \frac{e^{\xi_0}}{8N(N-2)} + \left(\e - \frac{e^{\xi_0}}{4N^2}\right)r^2 - \frac{e^{\xi_0}}{N^2(N^2 - 4)}r^{2 - N}.
\end{align*}
Combining the above estimates, we conclude that if $0 < \e \leq \e_1 := \min(1, \frac{e^{\xi_0}}{4N^2})$,
\begin{align*}
\Delta^{m-2} u(r) \leq -\beta + \frac{e^{\xi_0}}{2N} =: h(\beta) \quad \mbox{for any }\; r \in [0, \infty).
\end{align*}
This yields then for $\e \leq \e_1$, by Young's inequality,
\begin{align*}
u(r) \leq \e + \e\sum_{k=1}^{m-3} \frac{(-1)^k}{c_k}r^{2k} + h(\beta)\frac{r^{2m-4}}{c_{m-2}} \leq 2\e_1 + \Big[C_1 + h(\beta)\Big]\frac{r^{2m-4}}{c_{m-2}}, \quad \forall\; r > 0.
\end{align*}
As $\lim_{\beta\to\infty} h(\beta) = -\infty$, there exists $\beta_1$ large such that $u(r) \leq \ln\l_{N, m} - 2m\ln r$ in $(0, \infty)$ if $\beta \geq \beta_1$. This means that $u$ is stable for any $0 <\e \leq \e_1$ and $\beta \geq \beta_1$.
\qed

Our next result is inspired by \cite{DG}, where we construct some stable solutions to \eqref{1.1} by super-sub solution method.
\begin{thm}\label{thm1.5}
For any $m\geq 2$ and $N>2m$, let $P(x)$ be a polynomial verifying $$\lim_{|x|\to\infty} \frac{P(x)}{\ln|x|} = \infty \quad \mbox{and} \quad {\rm deg}(P)\leq 2m-2.$$ Then there exists $C_P \in \R$ such that for any $C \geq C_P$, we have a solution $u$ of \eqref{1.1} verifying
$$
- P(x) - C \leq u(x) \leq -P(x) - C + (1+|x|^2)^{m-\frac{N}{2}}\quad \mbox{in } \; \R^N.
$$
Consequently, there exists $\widetilde C_P \in \R$ such that the above solution $u$ is stable in $\mathbb{R}^N$ for any $C \geq \widetilde C_P$.
\end{thm}

\begin{rem}
We do not know if the assumption $\lim_{|x|\to\infty} \frac{P(x)}{\ln|x|} = \infty$ is equivalent or not to the apparently weaker condition $\lim_{|x|\to\infty} P(x) = \infty$.
\end{rem}

\noindent
Proof. We are looking for a solution $u$ of the form $u(x) = -P(x)-C+z(x)$ with
\begin{equation}\label{1.11}
(-\Delta)^m z(x) =e^{-P(x)-C +z(x)} \;\; \mbox{in }\;\R^N \quad \mbox{and}\quad z(x) = O(|x|^{2m-N}) \;\mbox{as } |x| \to\infty.
\end{equation}
Equivalently, we will resolve the following system:
\begin{equation}
\label{1.12}
\left\{
\begin{aligned}
&-\Delta z = (N-2m)(2m -2)v_1&&\text{in}~\mathbb{R}^N,\\
&-\Delta v_k = (N-2m+2k)(2m-2k-2)v_{k+1} &&\text{in}~\mathbb{R}^N, \; 1\leq k \leq m-2\\
&-\Delta v_{m-1} = d_m e^{-P(x)-C}e^z &&\text{in}~\mathbb{R}^N.
\end{aligned}
\right.
\end{equation}
Here
$$\frac{1}{d_m} = \prod_{i = 1}^{m-1} 2i(N-2i-2).$$

Set $W_j :=(1+|x|^2)^{j-\frac{N}{2}}$ for $j \in \Z$, the straightforward calculations yield that
$$
-\Delta W_j=(N-2j)(2j-2)W_{j-1}+(N-2j)(N-2j+2)W_{j-2} \quad \mbox{for any } j \in \Z.
$$
Therefore, for $2\leq j < \frac{N}{2}$, we have $-\D W_j \geq (N-2j)(2j-2)W_{j-1}$.

\medskip
Let $N > 2m$,
$$
Z(x):= W_m(x) > 0,\quad V_k := W_{m-k}(x) > 0\; \text{ for }\; 1\leq k \leq m-1.
$$
So $-\Delta Z\geq (N-2m)(2m-2)V_1$, $-\Delta V_k \geq (N-2m+2k)(2m-2k-2)V_{k+1}$ for $1\leq k\leq m-2$
and
$$
-\Delta V_{m-1} = N(N-2)W_{-1} = N(N-2)(1+|x|^2)^{-1-\frac{N}{2}}.
$$
Consider $$
f(x) := -P(x) + \frac{N+2}{2}\ln(1+|x|^2) + \ln d_m -\ln[N(N-2)]+(1+|x|^2)^{m-\frac{N}{2}},
$$
by our assumption on $P$ and $m < \frac{N}{2}$, readily $\max_{\R^N} f(x) = C_P < \infty$ exists. For any $C \geq C_P$, we have
$$
 -\Delta V_{m-1}\geq d_m e^{-P(x)-C_P}e^Z \geq d_m e^{-P(x)-C}e^Z \quad\text{in }~\mathbb{R}^N.
$$
In other words, $(Z, V_1,..., V_{m-1})$ is a super-solution in $\R^N$ to the system \eqref{1.12} for $C \geq C_P$.

\medskip
Since the system \eqref{1.12} is cooperative, $(0,0,...,0)$ and $(Z, V_1,..., V_{m-1})$ form a pair of ordered sub and super-solutions, we obtain the existence of a solution to \eqref{1.12}, hence a solution of \eqref{1.11}. Moreover, the solution $u$ satisfies $-P(x) - C \leq u(x) \leq -P(x)-C + Z(x)$ in $\R^N$.

\medskip
To ensure the stability of $u$, it's sufficient to choose $C$ such that
\begin{equation}\label{1.13}
 e^{u(x)} \leq e^{-P(x)-C + Z(x)} \leq e^{-P(x)- C + 1}\leq \frac{\gamma_{N,m}}{|x|^{2m}} \quad\text{in }~\mathbb{R}^N,
\end{equation}
 where $\gamma_{N,m} = \lambda_{N,m}$ in \eqref{1.3} if $m$ is odd and $\gamma_{N,m} = \mu_{N,m}$ given by \eqref{1.4} if $m$ is even. Let $
 g(x)=1-\ln\gamma_{N,m}-P(x)+2m\ln|x|$, clearly $C_P' = \max_{\R^N\backslash\{0\}} g(x)<\infty$ exists since
 $$
 \lim_{|x|\to 0} g(x) = \lim_{|x|\to\infty} g(x)=-\infty.
 $$
 Therefore, if we take $\widetilde C_p = \max(C_P, C_P')$, $u$ is a stable solution in $\mathbb{R}^N$ if $C \geq \widetilde C_P$.\qed

 \medskip
An immediate consequence of the above result is
 \begin{cor}\label{cor1.6}
Fro any $m\geq 3$ and $N>2m$, there exist non radial stable solutions to \eqref{1.1}. Moreover, when $m \geq 4$ is even, there are radial stable solutions on the borderline hypersurface of existence, i.e.~when $a_{m-1} = \Phi(a_k)$.
 \end{cor}

 \noindent
 Proof. Indeed, if $P$ is non radial in Theorem \ref{thm1.5}, the solution $u$ constructed is clearly non radial. On the other hand, if $P$ is radial, as our super and sub-solutions are radial, we can work in the subclass of radial functions to get a radial solution $u$. So for $m\geq 4$ even, if we consider polynomials $P(r) = \sum_{0\leq k\leq j}b_kr^{2k}$ with $b_j > 0$ and $1\leq j \leq m-2$, we obtain radial stable solutions $u$ satisfying $u(r) = o(r^{2m-2})$ at infinity. By \cite{FF}, such radial solutions must be on the borderline hypersurface $a_{m-1} = \Phi(a_k)$. \qed

\begin{rem}
For $m \geq 3$ odd, if we take $P(x) = P(r) = b_1r^2$ with $b_1 > 0$, the radial stable solutions
obtained verify that $(-\D)^ku(0) > 0$, i.e. sign$(a_k) = (-1)^k$, since otherwise $u(r) \leq -Cr^4$ at infinity, see \cite{FF}. The solutions obtained in the proof of Theorem \ref{thm1.2} are different, because they satisfy $\lim_{r\to\infty}\D^{m-1}u < 0$.
 \end{rem}

Our proof of Theorem \ref{thm1.4} is based on the following estimate.
\begin{lem}
\label{lem3.1}
Let $\xi$ be a radial function in $C^2(\mathbb{R}^N)$. Suppose that $\Delta \xi\geq r^\ell g(r)$ with $\ell>-1$ and $g$ nonincreasing in $r$, then
$$
\xi (r)\geq \xi(0)+\frac{r^{\ell+2}}{(N+\ell)(\ell+2)}g(r),\quad \forall~r\geq0.
$$
\end{lem}
In fact, we have
\begin{equation}\label{new7}
\xi'(r)\geq r^{1-N}\int_0^rg(s)s^{N-1} s^\ell ds
\geq r^{1-N}g(r)\int_0^rs^{N+\ell-1} ds
=\frac{r^{\ell+1}}{N+\ell} g(r).
\end{equation}
Integrating again, we get
$$
\xi(r)\geq\xi(0)+g(r)\frac{r^{\ell+2}}{(N+\ell)(\ell+2)}.
$$

\noindent
Proof of Theorem \ref{thm1.4}. Let $m$ be odd and $u$ be the solution to \eqref{1.2} with $a_k \leq 0$,
$1\leq k \leq m-1$.

\medskip
Let $w_k= \Delta^k u$, $k=1,...,m-1$. As $\Delta^{m-1} w_1 = -e^u < 0$ and $\Delta^k w_1(0) = a_{k+1} \leq 0$ for all $0
\leq k \leq m-2$, we get $w_1 \leq 0$ in $\R^N$, hence $u$ is decreasing in $r$. By Lemma \ref{lem3.1},
as $-\Delta w_{m-1} = e^u$,
$$
\begin{aligned}
-w_{m-1}(r) \geq -a_{m-1}(0)+ \frac{r^2}{2N}e^{u(r)}\geq\frac{r^2}{2N}e^{u(r)},
\end{aligned}
$$
so we have
$$
-\Delta w_{m-2}(r) = -w_{m-1}(r) \geq\frac{r^2}{2N}e^{u(r)}, \quad \forall~r>0.
$$
Applying again Lemma \ref{lem3.1}, we obtain
$$
-w_{m-2}(r) \geq -a_{m-2} +  \frac{r^4}{8N(N+2)}e^{u(r)} \geq  \frac{r^4}{8N(N+2)}e^{u(r)}.
$$
By induction, for all $1\leq k \leq m-1$,
$$
-w_{m-k}(r)\geq \frac{r^{2k}}{P_k(N)}e^{u(r)}\quad \mbox{for any } r>0,
$$
where
$$P_k(N)=2^k k!\prod_{\ell=0}^{k-1}(N+2\ell).$$
In particular, there holds
$$
-\Delta u(r) = -w_1(r) \geq \frac{r^{2m-2}}{P_{m-1}(N)}e^{u(r)}, \quad \forall~r>0.
$$
Using \eqref{new7}, we get
$$
-u'(r) \geq \frac{r^{2m-1}}{(N+2m-2)P_{m-1}(N)}e^{u(r)}, \quad \forall~r>0.
$$
Therefore
$$e^{-u(r)} \geq e^{-u(0)} + \int_0^r \frac{s^{2m-1}}{(N+2m-2)P_{m-1}(N)}ds \geq \frac{r^{2m}}{P_m(N)},$$
hence $$e^{u(r)}\leq \frac{P_m(N)}{r^{2m}} \quad \mbox{for any } r > 0.$$

As polynomial in $N$, deg$(P_m)=m$ while deg$(\lambda_{N,m})=2m$, so there exists $N_0$ such that for $N\geq N_0$, $P_m(N) \leq \lambda_{N,m}$, then $e^u\leq \frac{P_m(N)}{r^{2m}}\leq\frac{\lambda_{N,m}}{r^{2m}}$ i.e.~the solution $u$ is stable in $\mathbb{R}^N$.

\bigskip
\noindent

\end{document}